\documentclass[12pt]{article}
\usepackage[margin=0.8in]{geometry}
\usepackage{amssymb,eucal,xypic}
\usepackage[utf8]{inputenc}
\usepackage[english]{babel}
\usepackage{amsfonts}
\usepackage{amsmath}
\usepackage{amsthm}
\usepackage{mathtools}
\usepackage[all]{xy}
\usepackage{mathrsfs}
\usepackage{enumerate}
\usepackage{tikz-cd}
\usepackage{pst-node}
\usepackage{auto-pst-pdf}
\newtheorem{The}{Theorem}[section]
\newtheorem{Lem}[The]{Lemma}
\newtheorem{Prop}[The]{Proposition}
\newtheorem{Cor}[The]{Corollary}
\newtheorem{Rem}[The]{Remark}
\newtheorem{Obs}[The]{Observation}

\newtheorem{Question}[The]{Question}

\newtheorem{Def}[The]{Definition}

\newtheorem{Not}[The]{Notation}

\usepackage{xcolor}
\title{PhD}
\author{Erfan Soheil}
\date{September 2022}
\begin{document}

\begin{center}

{\Large\bf Properties of Critical Points of the Dinew-Popovici Energy Functional}

\end{center}

\begin{center}
{\large Erfan Soheil}

\end{center}


\begin{abstract}
Recently, Dinew and Popovici introduced and studied an energy functional $F$ acting on the metrics in the Aeppli cohomology class of a Hermitian-symplectic metric and showed that in dimension 3 its critical points (if any) are K\"ahler. In this article we further investigate the critical points of this functional in higher dimensions and under holomorphic deformations. We first prove that being a critical point for $F$ is a closed property under holomorphic deformations. We then show that the existence of a K\"ahler metric $\omega_k$ in the Aeppli cohomology class is an open property under holomorphic deformations. Furthermore, we consider the case when the $(2,\,0)$-torsion form $\rho_{\omega} ^{2,\,0}$ of $\omega$ is $\partial$-exact and prove that this property is closed under holomorphic deformations. Finally, we give an explicit formula for the differential of $F$ when the $(2,\,0)$-torsion form $\rho_{\omega} ^{2,\,0}$ is $\partial$-exact. 

\end{abstract}

\vspace{1ex}
\section{Introduction}\label{section:Introduction}
Let $X$ be a compact complex manifold with $\text{dim}_{\mathbb{C}} X=n$ and $\omega$ a Hermitian metric on $X$. This means that $\omega \in C^{\infty} _{1,\,1}(X,\mathbb{C})$ and $\omega > 0$. Let recall the following standard definitions.  
\begin{Def}\label{metrics}
\begin{itemize}
\item[(i)] $\omega$ is called K\"ahler if $d \omega=0$. We say that $X$ is a K\"ahler manifold if there exits a K\"ahler metric $\omega$ on $X$.
\item[(ii)] $\omega$ is called Hermitian-symplectic (H-s) if there exists $\rho^{2,\,0} \in C^{\infty}_{2,\,0}(X,\,\mathbb{C})$ such that 
\begin{equation}\label{eqn:H-S_condition}d(\rho^{2,\,0}  + \omega + \rho^{0,\,2}) = 0,\end{equation} where $\rho^{0,\,2}:=\overline{\rho^{2,\,0}}$. We denote $\Omega=\rho^{2,\,0} + \omega + \rho^{0,\,2}$ the \textbf{corresponding completion} of $\omega$. We say that $X$ is a Hermitian-symplectic manifold if there exists a Hermitian-symplectic metric $\omega$ on $X$.
\item[(iii)] $\omega$ is called SKT (pluriclosed) if $\partial\bar\partial \omega=0$. We say that $X$ is a SKT manifold if there exists a SKT metric $\omega$ on $X$.
\item[(iv)] $\omega$ is called balanced if $ d\omega^{n-1}=0$. We say that $X$ is a balanced manifold if there exists a balanced metric $\omega$.
\end{itemize}
\end{Def}

By a holomorphic family of compact complex manifolds we mean a proper holomorphic submersion $\pi: \mathcal{X}\rightarrow B$ between complex manifolds $\mathcal{X}$  and $ B$. This means that for every $t \in B$, $X_t=\pi^{-1}(t)$ is a compact complex submanifold of $\mathcal{X}$. From now on we denote $(X_t)_{t \in B}$ as a holomorphic family of compact complex manifolds instead of referring to $\pi: \mathcal{X}\rightarrow B$. If $B$ is simply connected then by Ehresmann's theorem (see \cite{Ehr73}
) all fibers $X_t$ are diffeomorphic. So $\mathcal{X}$ can be considered as a $C^{\infty}$ manifold $X$ equipped with a holomorphic family $(J_t)_{t \in B}$ of complex structures $((X,\, (J_t)_{t \in B}))$. From now on $B$ is an open ball containing the origin in $\mathbb{C}^m$. \\
One of the key theorems in deformation theory is the following statement by Kodaira and Spencer. The statement of the theorem is as follows.
\begin{The}\label{Thm: Kodaira Kahler open}(\cite{Kod3}, Theorem 15)
Suppose $(X_t)_{t \in B}$ is a holomorphic family of compact complex manifolds. If $X_0$ is a K\"ahler manifold, then any $X_t$ for all $t$ close enough to $0$ is again a K\"ahler manifold.
\end{The}
So the property of being a K\"ahler manifold is open under holomorphic deformations.
But the class of K\"ahler metrics is not the only class of metrics that is open under holomorphic deformations.\\

In \cite{pop10}, Popovici showed that the strongly Gauduchon property is open under holomorphic deformations as well. The notion of a strongly Gauduchon manifold was introduced by Popovici in \cite{Pop09}. Recall that $\omega$ is called strongly Gauduchon if $\partial\omega^{n-1}$ is $\bar\partial$-exact and we say that $X$ is said to be strongly Gauduchon manifold if there exists a strongly Gauduchon metric on $X$.\\

However, the openness property for an arbitrary class of metrics does not always hold. As a famous example, consider a holomorphic family of compact complex manifolds $(X_t)_{t \in B}$, and suppose $\omega_0$ is a balanced metric on $X_0$. \\
In \cite{Ale}, it is shown that the balanced property is not open under holomorphic deformations. Alessandrini and  Bassanelli pointed out the counter-example of the Iwasawa manifold endowed with the holomorphically parallelizable complex structure.\\  
Another example is that of class $\mathcal{C}$ manifolds. $X$ is of class $\mathcal{C}$ if it is bimeromorphic to a compact K\"ahler manifold. A counter-example was observed by Campana in \cite{cam91}.\\

Another class of manifolds that has drawn a lot of attention is the one of $\partial\bar\partial$-manifolds because they satisfy the Hodge decomposition and the Hodge symmetry. Recall that $X$ is called a $\partial\bar\partial$-manifold if and only if for every $d$-closed pure-type form $u$ on $X$ the following exactness properties are equivalent (the conclusion of $\partial\bar\partial$-lemma):
\begin{equation*}
u \: ~\text{is}~\: d\text{-exact}\Leftrightarrow \: u ~\: \text{is}~\:~ \partial\text{-exact}\Leftrightarrow ~\: u \:~ \text{is}\: ~\bar\partial\text{-exact}\Leftrightarrow ~\: u ~\: \text{is} ~\: \partial\bar\partial\text{-exact}.
\end{equation*} 

In  \cite{Wu} C.C. Wu proved that the $\partial\bar\partial$-property is open under holomorphic deformations.\\
In fact, if one considers a holomorphic family of compact complex manifolds $(X_t)_{t \in B}$ and supposes that the central fiber $X_0$ is a $\partial\bar\partial$-manifold then both the SKT and the  balanced properties become open under holomorphic deformations. However, in general, being an SKT manifold is not an open property under holomorphic deformations (see \cite{Fino}). 
\\

In a more general setting, we do not consider our manifolds to be $\partial\bar\partial$-manifolds. The main class of metrics that we discuss in this article is that of Hermitian-symplectic metrics. In dimension 2 any Hermitian-symplectic metric is K\"ahler (see \cite{street}) but in higher dimensions, the following question is still open.
\begin{Question}(\cite{street}, Question 1.7]) Do there exist non-Kähler Hermitian-symplectic complex
manifolds $X$ with $\text{dim}_{\mathbb{C}} X \geqslant 3$?
\end{Question}
Also in \cite{bell}, H. Bellitir, proved that the property of having a Hermitian-symplectic metric is open under holomorphic deformations.\\ In Definition \ref{metrics} (ii) $\Omega$ is not of type $(1,1)$ and $\rho^{2,\,0}$ is not unique. 
One can find a unique  $(2,\,0)$-form such that has the minimal $L^2 _{\omega}$-norm among such all forms, which we call the $(2,0)$-\textbf{torsion form} of $\omega$ and it is denoted by  $\rho_{\omega} ^{2,\,0}$. \\

 The main discussion of this article is based on \cite{pop20}, where Dinew and Popovici introduced the  \textbf{Dinew-Popovici energy functional}. Let $\omega_0$ be a fixed Hermitian-symplectic metric on $X$. They define ${\cal S}_{\{\omega_0\}}$  as follows
 $${\cal S}_{\{\omega_0\}}:=\{\omega_0 + \partial\bar{u_0} + \bar\partial u_0 \,\mid\,u_0 \in C^\infty_{1,\,0}(X,\,\mathbb{C}) ~~\text{such that}~~ \omega_0 + \partial\overline{u_0} + \bar\partial u_0>0\}.$$
The definition of Dinew-Popovici energy functional $F$ is given by
\begin{equation}\label{eqn: energy functional}
F : {\cal S}_{\{\omega_0\}} \to [0,\,+\infty), \hspace{3ex} F(\omega) = \int\limits_X|\rho_\omega^{2,\,0}|^2_\omega\,dV_\omega = ||\rho_\omega^{2,\,0}||^2_\omega,
\end{equation}
where $\omega \in {\cal S}_{\{\omega_0\}} $ and $\rho_\omega^{2,\,0}$ is the $(2,0)$-torsion form of $\omega$, while $|\,\,\,|_\omega$ is the pointwise norm and $||\,\,\,||_\omega$ is the $L^2$ norm induced by $\omega$. \\
When the dimension of $X$ is 3, the critical points for $F$ are exactly the K\"ahler metrics in the Aeppli cohomology class of $\omega_0$. In Theorem \ref{The: Main} we show that this property is open under holomorphic deformations in any dimension. In other words, we prove the following 
\begin{The}\label{The: Main}
Suppose $B $ is an open ball in $\mathbb{C}^m$ containing the origin and $(X_t)_{t \in B}$ is a holomorphic family of compact complex manifolds of complex dimension $n$ satisfying the followings conditions:
\begin{itemize}
\item[1)] for every $t \in B$, $X_t$ is equipped with a Hermitian-symplectic metric $\omega_t$ and the family $(\omega_t)_{t \in B}$ is a $C^{\infty}$-family of $(1,\,1)$-forms,
\item[2)] for  $t = 0$, $\omega_0$ is a K\"ahler metric on $X_0$. 
\end{itemize}
Then after possibly shrinking $B$ about $0$, there exists a family of $(1, \, 1)$-forms  $(\tilde{\omega}_t)_{t \in B}$ such that 
\begin{itemize}
\item[a)] $\tilde{\omega}_t \in \lbrace\omega_t\rbrace_{A}$, where $\lbrace\omega_t\rbrace_{A}$ is the Aeppli cohomology class of $\omega_t$,
\item[b)] $\tilde{\omega}_t$ is a K\"ahler metric on $X_t$ for every $t \in B$,
\item[c)] $\tilde{\omega}_0= \omega_0$,
\item[d)] $(\tilde{\omega}_t)_{t \in B}$ is a $C^{\infty}$ family of metrics. 
\end{itemize} 
\end{The}
By Theorem \ref{Thm: Kodaira Kahler open}, the open property for K\"ahler metrics is known. But the way that we constructed the $C^{\infty}$ family of K\"ahler metrics $(\tilde{\omega}_t)_{t \in B}$ is different. The new result of Theorem \ref{The: Main} is that we have constructed a K\"ahler metric in a specific Aeppli cohomology class.\\
In higher dimension, $\text{dim}_{\mathbb{C}} X> 0$ the following question is still open
\begin{Question}
When $\text{dim}_{\mathbb{C}}X > 3$, are the critical points of the Dinew-popovici energy functional $F:{\cal S}_{\{\omega_0\}} \to [0,\,+\infty)$ exactly the K\"ahler metrics in the Aeppli cohomology class of $\omega_0$? 
\end{Question}
We give a partial answer to this question in Proposition \ref{intro :prop: critial d-exact} and Corollary \ref{critical dbar exact}. Precisely we show that
\begin{Prop}\label{intro :prop: critial d-exact}
Suppose that $(X,\omega_0)$ is a compact complex Hermitian-symplectic manifold of dimension $n$. Fix an $\omega \in {\cal S}_{\{ \omega_0 \}}$. If $\rho_\omega^{2,\,0}=\partial \xi$, for some $(1,0)$-form $\xi$, then the differential at $\omega$  of the Dinew-Popovici energy functional $F$ defined in equation (\ref{eqn: energy functional}) evaluated on $\gamma=\bar\partial\xi+\partial\bar{\xi}$ is
\begin{equation}\label{der F wrt xi}
d_{\omega}F(\gamma)= 2 \Vert \rho^{2,0} _{\omega} \Vert^2+ 2\mbox{Re} \int\limits_X  \bar\partial \xi \wedge \rho^{2,0} \wedge  \overline{\rho_\omega^{2,\,0}} \wedge \omega_{n-3}.
\end{equation} 
\end{Prop}
From this we get the following
\begin{Cor}\label{critical dbar exact}
Under the assumptions of Proposition \ref{intro :prop: critial d-exact} if
\begin{itemize}
\item[(i)] $\omega$ is a critical point for $F$, and
\item[(ii)]the $(2,\,0)$-torsion form $\rho_\omega^{2,\,0}=\partial\xi$ such that $\bar\partial \xi$ is weakly semi-positive,
\end{itemize}
then $\omega$ is a K\"ahler metric on $X$.
\end{Cor}
Moreover, in Proposition \ref{prop: close critical F_t}, we prove that the property of being a critical point for $F$ is closed under holomorphic deformations. Precisely, we prove the following proposition.
\begin{Prop}\label{prop: close critical F_t}
Suppose  $(X_t)_{t \in B}$ is a holomorphic family of compact complex manifolds, $(\omega_t)_{t \in B}$ is a $C^{\infty}$ family of Hermitian-symplectic metrics on $(X_t)_{t \in B}$ and $(F_t)_{t \in B}$ is the associated family of Dinew-Popovici linear functionals $F_t:{\cal S}_{\{\omega_t \}}\rightarrow [0,\, \infty]$ (see section \ref{section: back energy functioal}). If after possibly shrinking $B$ about $0$,
\begin{itemize}
\item[(1)] for every $t \in B\backslash \lbrace 0\rbrace$, $\omega_t$ is a critical point in $F_t$, 
\item[(2)] for every $t \in B$, $h_{BC,\,t} = h_{BC,\,0} $, where $h_{BC,\,t} $ is the dimension of $\ker H_{BC} ^{0,\,2} (X_t,\, \mathbb{C})$,
\end{itemize}
Then $\omega_0$ is a critical point for $F_0$. 
\end{Prop}
In the above statement, $ H_{BC} ^{0,\,2} (X_t,\, \mathbb{C})$ is the Both-Chern cohomology group of bidegree $(0,\,2)$ of $X_t$ (see Definition \ref{cohomologies}). In section \ref{section:Preliminaries} we first give the definitions and tools to state the main results and in section \ref{Result} we state our new results and prove them.
\section*{Acknowledgement}
I wish to thank Dan Popovici, my research supervisor, for his patience, support, and significant and essential helpful comments and criticisms of this study effort. This article would not be completed without his help..
\section{Preliminaries}\label{section:Preliminaries}
In this section, we recall the required definitions, lemmas, and propositions that  will be frequently used in section \ref{Result}. \\
Throughout this section,  $X$ is a compact complex manifold of dimension $n$ equipped with a Hermitian metric $\omega$. This means that $\omega$ is a $C^{\infty}$ positive definite $(1,\,1)$-form on $X$.
\subsection{General background on complex geometry}\label{sec back1}
This subsection contains some standard and well-known definitions and results  in complex geometry. The reader is referred to  \cite{Dem}, \cite{Huy}, and \cite{Voi} for further details. \\
First, we recall four different notions of positivity for differential forms.
Let $V$ be a complex vector space of dimension $n$ and $(z_1,\dots,z_n)$ be a coordinate on $V$. We denote the corresponding basis of $V$ by $(\partial/\partial z_1,\dots,\partial/\partial z_n)$ and its dual basis in $V^*$ by $(dz_1,\dots,dz_n)$. consider the exterior algebra
$$
\Lambda V^*_\mathbb C=\bigoplus\Lambda^{p,q}V^*,\quad \Lambda^{p,q}V^*=\Lambda^p V^*\otimes\Lambda^q\overline{V^*}.
$$
Since $V$ is a complex vector space, it has a canonical orientation, given by the $(n, n)$-form
$$
\tau(z)=idz\wedge d\bar z_1\wedge\cdots\wedge idz_n\wedge d\bar z_n=2^n dx_1\wedge dy_1\wedge\cdots\wedge dx_n\wedge dy_n,
$$
where $z_j=x_j+iy_j$. In fact, if $(w_1,\dots,w_n)$ are other the coordinates, we find
$$
dw_1\wedge\cdots\wedge d w_n=\det(\partial w_j/\partial z_k)dz_1\wedge\cdots\wedge dz_n,
$$
$$
\tau(w)=|\det(\partial w_j/\partial z_k)|^2\tau(z).
$$
So one can define the notion of positivity as independent of local coordinates. 
\begin{Def}\label{positivty}
\begin{itemize}
\item[(1)]
A $(q,q)$-form $v\in\Lambda^{q,q}V^*$ is said to be \textbf{strongly semi-positive} (resp. \textbf{strongly strictly positive}) if $v$ is a convex combination
$$
v=\sum\gamma_si\alpha_{s,1}\wedge\bar\alpha_{s,1}\wedge\cdots\wedge i\alpha_{s,q}\wedge\bar\alpha_{s,q}
$$
where $\alpha_{j,s}\in V^*$ and $\gamma_s\ge 0$ (resp. $\gamma_s > 0$ ).

\item[(2)]
A $(p,p)$-form $u\in\Lambda^{p,p}V^*$ is said to be  \textbf{weakly semi-positive} (resp.  \textbf{weakly strictly positive}) if for all $\alpha_j\in V^*$, $1\le j\le q=n-p$, then
$$
u\wedge i\alpha_1\wedge\bar\alpha_1\wedge\cdots\wedge i\alpha_q\wedge\bar\alpha_q \geq 0 ~~ (\text{resp.} ~~ u\wedge i\alpha_1\wedge\bar\alpha_1\wedge\cdots\wedge i\alpha_q\wedge\bar\alpha_q > 0)
$$
\end{itemize}
\end{Def}
\begin{Rem}\label{rem: positivity notions} Locally any Hermitian metric $\omega$ is a strongly strictly positive  $(1,\,1)$-form and $\omega$ has the following representation
$$\omega= \sum i dz_i \wedge d\bar{z}_i.$$ 
\end{Rem}
Fortunately, the concepts of weakly semi-positive (resp. weakly strictly positive) and strongly semi-positive (resp. strongly strictly positive) coincide in bidegree  $(1,1)$ and $(n-1,n-1)$. 
\begin{Prop}(\cite{Dem}, Chapter III, Proposition 1.11)\label{prod of sp} If $u_1,\dots,u_s$ are strongly semi-positive (resp. strongly strictly positive) forms, then $u_1\wedge\cdots\wedge u_s$ is strongly semi-positive (resp. strongly strictly positive) form.
\end{Prop}
For simplicity we recall the following notation.
\begin{Not}\label{Not: metric}
For any $k \in \mathbb{N}$, 
$$\omega_{k}=\frac{\omega^k}{k!}.$$
\end{Not}
It is obvious that $\partial\omega_k=\partial\omega \wedge \omega_{k-1}$ and $\bar\partial\omega_k=\bar\partial\omega \wedge \omega_{k-1}$. Also, it is a well-known fact that for a Hermitian metric $\omega$ we have
\begin{equation}\label{star op on omega }
\star_{\omega} \omega_k= \omega_{n-k},  ~~~~~~ \forall~ k \in \{1,\ \dots ,\ n\},
\end{equation}
where $\star_{\omega} $ is the Hodge star operator induced by $\omega$. The following proposition plays an important role in our discussion later.  
\begin{Prop}\label{primi formula for fomrs}(\cite{Voi}, Proposition 6.29)
If $u \in C^{\infty}_{p,\,q}(X,\,\mathbb{C})$ is primitive then 
\begin{equation}\label{eqn:primi formula for fomrs}
\star u = (-1)^{\frac{(p+q)^2 + p+q}{2}}i^{p-q} \omega_{n-q-p}\wedge u. 
\end{equation}
\end{Prop}
Recall that a $(p,\,q)$-form $u$ is primitive if $L_{\omega} ^{\star}(u)=0$, where $L_{\omega} ^{\star}$ is the adjoint of the Lefschetz operator $L_{\omega} (u)= \omega \wedge u$.\\
 Now we mention four equations which one can easily imply by equation (\ref{eqn:H-S_condition}).  
\begin{Obs}\label{obs: h-s dclosed}
If $\omega$ is a Hermitian-symplectic metric then
\begin{itemize}
\item[(i)] $\partial\omega = -\bar\partial\rho^{2,\,0} _\omega ~~ \text{and}~~ \bar\partial \omega= -\partial \rho^{0,\,2} _\omega$.
\item[(ii)] $\partial\rho^{2,\,0} _\omega =0 ~~\text{and}~~ \bar\partial \rho^{0,\,2} _\omega=0 $.
\item[(iii)] $\omega$ is K\"ahler if and only if $\rho^{2,\,0} _\omega =0$.
\item[(iv)]$\partial\bar\partial \omega=0$.
\end{itemize}
\end{Obs}
Note that (i) and (ii) imply that if $\omega$ is a Hermitian-symplectic metric then $\partial\omega$ and $\bar\partial\omega$ is $d$-closed.\\
In order to define suitable cohomology groups for Hermitian-symplectic and SKT metrics, we recall the definitions of the Bott-Chern cohomology and the Aeppli cohomology groups. 
\begin{Def}\label{cohomologies}
For every $p,q \in \lbrace  1,\,\ldots \, n \rbrace$ one defines:
\begin{itemize}
\item[(i)] the Bott-Chern cohomology group of bidegree (or type)  $(p,\, q)$ of $X$ as
\begin{equation}\label{bott-chern eqn}
H^{p,\,q} _{BC}(X ,\, \mathbb{C})= \frac{\ker \partial \cap \ker \bar\partial}{ \operatorname{Im}(\partial\bar\partial)},
\end{equation}
\item[(ii)]the Aeppli cohomology group of bidegree (or type) (p, q) of X as
\begin{equation}\label{aeppli eqn}
H^{p,\,q} _{A}(X ,\, \mathbb{C})=\frac{\ker(\partial\bar\partial)}{\operatorname{Im}\partial +\operatorname{Im}\bar\partial},
\end{equation}
\end{itemize}
where all the kernels and images are considered as $\mathbb{C}$-vector subspaces of $C^{\infty}_{p,\,q}(X,\,\mathbb{C})$ according to the case.
\end{Def}
From definitions \ref{metrics} and \ref{cohomologies} one can see if $\omega$ is a Hermitian-symplectic (resp. SKT) metric then the Aeppli (resp. Bott-Chern) cohomology class of $\omega$, which will be denoted by $\lbrace \omega\rbrace_{A}$ (respectively $\lbrace\omega\rbrace_{BC}$), is well-defined. \\
In the following definition, we recall formal definitions of two elliptic self-adjoint operators and mention the Hodge decompositions for $C^{\infty}_{p,\,q}(X,\,\mathbb{C})$ of these operators.  
\begin{Def}
Fix $p,\,q \in \lbrace  1,\,\ldots \, n \rbrace $ then
\begin{itemize}
\item[(i)] The Bott-Chern Laplacian operator $\Delta^{p, \, q}_{BC} : C^{\infty}_{p,\,q}(X,\,\mathbb{C}) \rightarrow C^{\infty}_{p,\,q}(X,\,\mathbb{C})$ is defined as follows
\begin{equation}\label{Delta-bt}
\Delta^{p, \, q}_{BC}:=\partial^{\star}\partial + \bar\partial^{\star}\bar\partial + (\partial\bar\partial)^{\star}(\partial\bar\partial) + (\partial\bar\partial)(\partial\bar\partial)^{\star} + (\partial^{\star}\bar\partial)^{\star}(\partial^{\star}\bar\partial) + (\partial^{\star}\bar\partial)(\partial^{\star}\bar\partial)^{\star},
\end{equation}
\item[(ii)] The Dolbeault Laplacian operator $\Delta^{p, \, q}_{\bar\partial} : C^{\infty}_{p,\,q}(X,\,\mathbb{C}) \rightarrow C^{\infty}_{p,\,q}(X,\,\mathbb{C})$ is defined as follows
\begin{equation}\label{Delta Dolbeault}
\Delta^{p, \, q}_{\bar\partial}:= \bar\partial\bar{\partial}^{\star}+\bar{\partial}^{\star}\bar\partial.
\end{equation}
\end{itemize}
\end{Def}
It is worth mentioning that by the definition of the  Bott-Chern Laplacian operators, it is a real self-adjoint operator but the conjugate Dolbeault 
operator is not. \\
For each of the above operators, we have the following $L^2 _{\omega}$ two-space orthogonal decomposition for $ C^{\infty}_{p,\,q}(X,\,\mathbb{C})$
\begin{itemize}
\item[(i)]\begin{equation}\label{ bott-cherndecompositions}
 C^{\infty}_{p,\,q}(X,\,\mathbb{C})= \ker\Delta^{p, \,q} _{\text{BC}}\oplus \operatorname{Im} \Delta^{p, \,q} _{\text{BC}}, 
\end{equation}
\item[(ii)]\begin{equation}\label{dolb decompositions}
 C^{\infty}_{p,\,q}(X,\,\mathbb{C})= \ker\Delta^{p, \,q} _{\partial}\oplus \operatorname{Im} \Delta^{p, \,q} _{\partial}.
\end{equation}
\end{itemize}


\subsection{Background on deformation of complex structures}\label{section : back deformation}
This subsection is a summary of some basic  definitions and results on  deformation of complex structures of compact complex manifolds. Our main references for this part are \cite{Kod86} and \cite{Kod3}. Also there are series of papers published by D. Popovici \cite{pop20}, \cite{Pop}, \cite{Pop19}, \cite{Pop14}, \cite{Pop09} and \cite{PSU20} which play a crucial role in this article, so we recall some lemmas and propositions from them. \\ 

We recall the green operator of a self-adjoint elliptic operator. For every fixed $p,\,q \in \lbrace 1,\,\ldots,\,n\rbrace$ suppose $E$ is a self-adjoint elliptic operator on $ C^{\infty}_{p,\,q}(X,\,\mathbb{C})$, since $X$ is a compact manifold $\ker E$ is a finite-dimensional complex vector space. We denote by $F_{E}: C^{\infty}_{p,\,q}(X,\,\mathbb{C}) \rightarrow \ker E$ the $L^2 _{\omega}$ orthogonal projection. One can define  the \textbf{Green operator} of $E$, $E^{-1} : C^{\infty}_{p,\,q}(X,\,\mathbb{C}) \longrightarrow \operatorname{Im} E$, such that 
\begin{equation}\label{green property}
E^{-1} E(\gamma)= E E^{-1} (\gamma)= \gamma - F_{E}(\gamma), ~~~~~ \gamma \in C^{\infty} _{p,\,q}(X,\,\mathbb{C}).
\end{equation}
If we restrict $E$ to  $\operatorname{Im} E$ then $E$ is a bijection and so $E^{-1} :  \operatorname{Im} E  \rightarrow \operatorname{Im} E$ is the inverse of this restriction. In particular one can define, $F_{\text{BC}}$, 
 $F_{\bar{\partial}}$, $\Delta^{-1} _{\text{BC}}$, 
and $\Delta^{-1} _{\bar{\partial}}$. \\
Thanks to \cite{pop20} we have all tools to give the explicit formula for the $(2,\,0)$-torsion form $\rho^{2,\,0} _{\omega}$ for any Hermitian-symplectic metric $\omega$. 
\begin{Lem}(\cite{pop20}, Lemma and Definition 3.1)\label{torsion fromula1}
Suppose $\omega$ is a  Hermitian-symplectic metric on $X$ and $\rho_\omega^{2,\,0}$ is the  $(2,\,0)$-torsion form of $\omega$. Then 
\begin{equation}\label{eqn:trosion}
\rho_\omega^{2,\,0} = -\Delta_{BC}^{-1}[\bar\partial^\star\partial\omega + \bar\partial^\star\partial\partial^\star\partial\omega].
\end{equation}
\end{Lem}
Notice that in (\ref{eqn:trosion}) by $\Delta_{BC}^{-1}$ we mean $(\Delta_{BC} ^{2,\,0})^{-1}: C^{\infty}_{2,\,0}(X,\,\mathbb{C}) \longrightarrow \operatorname{Im} \Delta_{BC} $. Also since $\Delta_{BC}$ is a real operator (so $\Delta_{BC}^{-1}$) we thus have
\begin{equation}\label{ean: conj torsion}
\rho_\omega^{0,\,2} = -\Delta_{BC}^{-1}[\partial^\star\bar\partial\omega + \partial^\star\bar\partial\bar\partial^\star\bar\partial\omega].
\end{equation}
The following theorem gives us a criteria to determine whether these families are $C^{\infty}$ family of linear operators. This is the main key to proving Theorem \ref{The: Main}.
\begin{The}\label{Kodaira}(\cite{Kod3}) \textit{Kodaira-Spencer fundamental theorem.}
\begin{itemize}
\item[(i)]If the dim $\ker \Delta_{\text{BC},\, t}: C^{\infty}_{p,\,q}(X_t,\,\mathbb{C}) \rightarrow C^{\infty}_{p,\,q}(X_t,\,\mathbb{C}) $ (resp. 
dim $\ker \Delta_{\bar\partial,\, t}$ ) is independent of $t \in B$, then the family $(F_{\text{BC},\,t})_{t \in B}$ (resp. 
 $(F_{\bar\partial,\,t})_{t \in B}$) is a $C^{\infty}$ family of linear operators.
\item[(ii)]If the dim $\ker \Delta_{\text{BC},\, t}: C^{\infty}_{p,\,q}(X_t,\,\mathbb{C}) \rightarrow C^{\infty}_{p,\,q}(X_t,\,\mathbb{C}) $ (resp. 
dim $\ker \Delta_{\bar\partial,\, t}$ ) is independent of $t \in B$, then the family $(\Delta^{-1}_{\text{BC},\,t})_{t \in B}$ (resp. 
$(\Delta^{-1}_{\bar\partial,\,t})_{t \in B}$) is a $C^{\infty}$ family of linear operators.
\end{itemize}
\end{The}
\subsection{Background on the Dinew-Popovici energy functional }\label{section: back energy functioal}
This subsection is devoted to some definitions and results  based on \cite{pop20}. Throughout this section $\omega_0$ is a Hermitian-symplectic metric on $X$ which we consider as our background metric. The main goal of this subsection is to find the explicit formula for differential of Dinew-Popovici energy functional $F$ in equation (\ref{eqn: energy functional}) at $\omega$, where $\omega \in {\cal S}_{\{\omega_0\}}$ is a fixed Hermitian-symplectic metric.    
  
\begin{Prop}\label{Prop:F-tilde_properties}(\cite{pop20}, Proposition 3.5) The differential at $\omega$ for $F$ is given by the formula: 
\begin{equation}\label{eqn:differential_F}
(d_{\omega}F)(\gamma) = -2\,\mbox{Re}\,\langle\langle u,\,\bar\partial^{\star}\omega\rangle\rangle_\omega + 2\,\mbox{Re}\,\int\limits_X u\wedge\rho_\omega^{2,\,0}\wedge\overline{\rho_\omega^{2,\,0}}\wedge\bar\partial\omega_{n-3}
\end{equation}
for every $(1,1)$-form $\gamma=\bar{\partial}u+\partial\bar u$.
\end{Prop}
In dimension $3$, the term $2\,\mbox{Re}\,\int\limits_X u\wedge\rho_\omega^{2,\,0}\wedge\overline{\rho_\omega^{2,\,0}}\wedge\bar\partial\omega_{n-3}$ vanishes and we have $(d_{\omega}F)(\gamma) = -2\,\mbox{Re}\,\langle\langle u,\,\bar\partial^{\star}\omega\rangle\rangle_\omega $. If $\omega$ is a critical point for $F$ and $u=\bar\partial^{\star}\omega$ then one can see that $\bar\partial^{\star}\omega=0$. This means that $\omega$ is a balanced metric, on the other hand $\omega$ is SKT. So We conclude that $\omega$ is K\"ahler.
Therefore in dimension 3 the critical points of $F$ are exactly K\"ahler metrics in the Aeppli cohomology class of $\omega_0$.\\
Now consider a holomorphic family of compact complex manifolds $(X_t)_{t \in B}$ and  $(\omega_t)_{t \in B}$ is a $C^{\infty}$ family of Hermitian-symplectic metrics on $(X_t)_{t \in B}$. This means that for every $t \in B$, $\omega_t$ is a Hermitian-symplectic metric on $X_t$. Therefore like (\ref{eqn: energy functional}) one can define a family of Dinew-Popovici linear functionals $(F_t)_{t \in B}$, $F_t : {\cal S}_{\{\omega_t\}} \rightarrow [0,\, \infty)$, as follows
\begin{equation}\label{eqn:F_energy-functional_H-S t}
F_t : {\cal S}_{\{\omega_t\}} \to [0,\,+\infty), \hspace{3ex} F_t(\bar\omega_t ) = \int\limits_{X_t}|\rho_{\bar\omega_t} ^{2,\,0}|^2_{\bar\omega_t}\,dV_{\bar\omega_t} = ||\rho_{\bar\omega_t} ^{2,\,0}||^2_{\bar\omega_t},  
\end{equation}
where like the Equation (\ref{eqn: energy functional}) $\bar\omega_t \in{\cal S}_{\{\omega_t\}}$ and $\rho_{\bar\omega_t}^{2,\,0}$ is the $(2,\,0)$-torsion form of $\bar\omega_t$, while $|\,\,\,|_{\bar\omega_t}$ is the pointwise norm and $||\,\,\,||_{\bar\omega_t }$ is the $L^2$ norm induced by $\bar\omega_t $.\\
Henceforth if we fix any Hermitian symplectic $\bar\omega_t \in {\cal S}_{\{\omega_t\}}$ then for every $t \in B$ one can define the differential at $\bar\omega_t$ of $F_t$ exactly like Proposition \ref{Prop:F-tilde_properties}.
\section{Results}\label{Result}
This  section is devoted to our new results based on \cite{pop20}. We give a proof for Theorem \ref{The: Main}.  
This theorem shows that if a compact complex manifold $X$ admits a Hermitian-symplectic metric $\omega_0$, then the existence of a K\"ahler metric $\tilde{\omega}_0 $ in the Aeplli cohomology class $\omega_0$ is an open property  under holomorphic deformations. 
Before we present the proof of Theorem \ref{The: Main}, we mention three theorems which play a crucial role in our proof. 
\begin{The}(\cite{Pop}, Theorem 4.1)\label{The:min-sol-ddbar} Fix a compact Hermitian manifold $(X,\,\omega)$. For any $C^{\infty}$ $(p,\,q)$-form $v\in\mbox{Im}\,(\partial\bar\partial)$, the (unique) minimal $L^2$-norm solution of the equation

\begin{equation}\label{eqn:ddbar-eq}\partial\bar\partial u=v\end{equation}

\noindent is given by the formula

\begin{equation}\label{eqn:ddbar-eq-solformula1}u = (\partial\bar\partial)^{\star}\Delta_{BC}^{-1}v,\end{equation}

\noindent where $\Delta_{BC}^{-1}$ is the Green operator of the Bott-Chern Laplacian $\Delta_{BC}$ induced by $\omega$.

\end{The}
\begin{The}\label{Ballanced opennes}(\cite{Wu}, Theorem 5.12)
Let $(X_t)_{t \in B}$ be a holomorphic family of compact complex manifolds of complex dimension $n$. If the central fiber $X_0$ is a $\partial \bar\partial$-manifold, then after possibly shrinking $B$ about $0$, $X_t$ is a $\partial \bar\partial$-manifold  for all $t \in B$. 
\end{The}
\begin{The}(\cite{deli}, Section 6)\label{kahler to ddbar }
Every compact K\"ahler manifold is a $\partial\bar\partial$-manifold.
\end{The}
\textbf{Proof of Theorem \ref{The: Main}}.
Since $\omega_0$ is a  K\"ahler metric on $X_0$ by Theorem \ref{kahler to ddbar }, $X_0$ is a $\partial\bar\partial$-manifold, therefore by Theorem \ref{Ballanced opennes} after possibly shrinking $B$ about $0$ one can assume that $X_t$ is a $\partial\bar\partial$-manifold for every $t \in B$. Let fix a $t \in B$, $\omega_t$ is a Hermitian-symplectic metric on $X_t$ then by Observation \ref{obs: h-s dclosed} in Section \ref{sec back1} one implies that  $\partial_t\omega_t$ is $d$-closed and $\partial_t$-exact. Since $X_t$ is a $\partial\bar\partial$-manifold, $\omega_t$ is $\partial_t\bar\partial_t$-exact. So the following equation 
\begin{equation}\label{exis}
-\partial_t\bar\partial_t u_t=\partial_t\omega_t.
\end{equation}
has at least one solution, $u_t$, for $t \in  B$. By Theorem \ref{The:min-sol-ddbar} we are able to choose the minimal $L^2$-norm solution with respect to $\omega_t$ among all such $u_t$. The minimal  $L^2 _{\omega_t}$-norm solution of equation (\ref{exis}) is given by
\begin{equation}\label{minkahler}
u_t ^{min}= -(\partial_t \bar\partial_t)^{\star}\Delta_{BC, t}^{-1} (\partial_t \omega_t),
\end{equation} 
where $\Delta_{BC,t}^{-1}$ is the Green operator of the Bott-Chern Laplacian $\Delta_{BC,\, t}$ induced by $\omega_t$, mentioned in Section \ref{section : back deformation}. Now we define,
 \begin{equation}\label{eqn: kahler_t}
 \tilde{\omega}_t= \omega_t+ \partial_t \overline{u_t ^{min}} + \bar\partial_t u_t ^{min},
\end{equation}
for all $t \in B$.\\
By the construction of $ \tilde{\omega}_t$, one can see that 
$$\partial_t\bar\partial_t \tilde{\omega}_t = \partial_t\bar\partial_t(\omega_t+ \partial_t \overline{u_t ^{min}} + \bar\partial_t u_t ^{min})=  \partial_t\bar\partial_t\omega_t =0.$$
Therefore $\lbrace \tilde{\omega}_t \rbrace_A$ is well-defined and by the definition of the Aeppli cohomology group, adding $\partial_t \overline{u_t ^{min}}$ and $ \bar\partial_t u_t ^{min}$ to $\omega_t$ does not change the Aeppli cohomology class of $\omega_t$. Hence $\tilde{\omega}_t \in \lbrace \omega_t \rbrace_{A}$, this proves (a). \\ Also for every $t \in B$, $\tilde{\omega}_t$ is $d$-closed because
\begin{equation}\label{d-closed tilde omega}
d\tilde{\omega}_t= d(\omega_t+ \partial_t \overline{u_t ^{min}} + \bar\partial_t u_t ^{min})=\partial_t\omega_t+\bar\partial_t\omega_t+\bar\partial_t\partial_t \overline{u_t ^{min}}+\partial_t\bar\partial_t u_t ^{min}.
\end{equation}
Equation (\ref{minkahler}) implies that $\partial\bar\partial_t u_t ^{min}=- \partial_t\omega_t$, put this in the equation (\ref{d-closed tilde omega}) one can see that $\tilde{\omega}_t$ is d-closed. On the other hand, the strict positivity of $\omega_0$ implies strict positivity of  $ \tilde{\omega}_t$  for all $t \in B$ sufficiently close to $0$, henceforth $\tilde{\omega}_t$ is a strictly positive $d$-closed $(1,\,1)$-form on $X_t$, this proves (b). \\ 
So we have a family of K\"ahler metrics $(\tilde{\omega}_t)_{t \in B}$ on $(X_t)_{t \in B}$. At $t=0$, there are two K\"ahler metrics on $X_0$. One of them is $\omega_0$,  which is given by assumption (2) and the other one is $\tilde{\omega}_0$ by our construction. Since $\omega_0$ is a K\"ahler metric on $X_0$, $\partial_0\omega_0=0$. Hence
$$u_0 ^{min}= -(\partial_0 \bar\partial_0)^{\star}\Delta_{BC,0}^{-1}(\partial_0 \omega_0)=0. $$
So,
$$\tilde{\omega}_0= \omega_0+ \partial_0 \overline{u_0 ^{min}} + \bar\partial_0 u_0 ^{min}= \omega_0. $$
This means that these two metrics coincide on $X_0$ which proves (c).\\
For every $t \in B$ we denote by $h_{\text{BC,\,t}}(X_t)$ the dimension of $\ker\Delta_{\text{BC},\, t}$ ($\Delta_{\text{BC},\, t}: C^{\infty}_{2,\,1}(X_t,\,\mathbb{C}) \rightarrow C^{\infty}_{2,\,1}(X_t,\,\mathbb{C}))$. Since $(X_t)_{t \in B}$, after possibly shrinking $B$ about $0$, is a holomorphic family of compact complex $\partial\bar \partial$-manifolds, by Theorem 5.12 in \cite{Wu}, $h_{\text{BC},\,t}(X_t)=h_{\text{BC},\,0}(X_0)$  for every $t \in B$. By Theorem \ref{Kodaira} (ii)  the family of linear operators $(\Delta^{-1} _{BC,\,t})_{t \in B}$ acting on $(2,\,1)$-forms is a $C^{\infty}$ family of linear operators, therefore the family $(u_t ^{min})_{t \in B}$ is a $C^{\infty}$ family of $(1, \, 0)$-forms, and since $(\omega_t)_{t \in B}$ is a $C^{\infty}$ family of metrics one can say $(\tilde{\omega}_t)$ is a $C^{\infty}$ family of metrics, this proves (d).\hfill $\Box$\\
We saw that in dimension $3$, K\"ahler metrics are the critical points for the Dinew-Popovici energy functional $F$. So as a consequence of Theorem \ref{The: Main}  one can get the following corollary.

\begin{Cor}\label{crit kahler}
Suppose  $(X_t)_{t \in B}$ is a holomorphic family of compact complex manifolds of dimension $3$, $(\omega_t)_{t \in B}$ is a $C^{\infty}$ family of Hermitian-symplectic metrics on $(X_t)_{t \in B}$, and $(F_t)_{t \in B}$ is a family of Dinew-Popovici linear functionals mentioned in section \ref{section: back energy functioal}. If  for $t = 0$, $\omega_0$ is a critical point of $F_0$, then after possibly shrinking $B$ about $0$ there exists a $C^{\infty}$ family of K\"ahler metrics $(\tilde{\omega}_t)_{t \in B}$ such that for every $t \in B$, $\tilde{\omega}_t \in {\cal S}_{\{\omega_t \}} $ and $\tilde{\omega}_t$ is a critical point of $F_t$ and $\tilde{\omega}_0=\omega_0$.
\end{Cor}
\textbf{Proof.} The existence of a $C^{\infty}$ family of K\"ahler metrics $(\tilde{\omega}_t)_{t \in B}$ such that for every $t \in B$ each $\tilde{\omega}_t \in {\cal S}_{\{\omega_t \}}$ and $\tilde{\omega}_0=\omega_0$ come directly from Theorem \ref{The: Main} and since the dimension of each $X_t$ is 3, $\tilde{\omega}_t$ being a K\"ahler for each $t \in B$ implies that $\tilde{\omega}_t$ is a critical point of $F_t$ . \hfill $\Box$ \\

By Corollary 4.2 of \cite{pop20}, in dimension $3$ if $\omega$ is a Hermitian-symplectic and the given Aeppli class $\lbrace\omega\rbrace_{A} $ contains a K\"ahler metric $\omega_k$, then its $(0,\,2)$-torsion form $\rho_{\omega} ^{0,\,2}$ is $\bar\partial$-exact. Therefore by Theorem \ref{The: Main} if the given Aeppli class $\lbrace\omega\rbrace_{A} $ contains a K\"ahler metric $\omega_k$, then the $\bar\partial$-exactness for $\rho_{\omega} ^{0,\,2}$  is an open property under holomorphic deformations.\\
So it is natural to investigate the openness and the closedness properties of the $(0,\,2)$-torsion form $\rho_{\omega} ^{0,\,2}$ in  higher dimensions.\\
In the following proposition we show that for a  Hermitian-symplectic metric $\omega$, the $\bar\partial$-exactness for the $(0,\,2)$-torsion form $\rho_{\omega} ^{0,\,2}$ is a closed property under small holomorphic deformations in any dimension.\\
First, we fix some notations for next proposition. For every $t \in B$ let $h_{BC,\, t} =\text{dim} \ker \Delta_{BC,\, t} ^{0,2}$ and $h_{\bar\partial,\, t} =\text{dim} \ker \Delta^{\ 0,2} _{\bar\partial,\,t}$
\begin{Prop}\label{prop: close torsion dbar exact}
Suppose that $(X_t)_{t \in B}$ is a holomorphic family of compact complex manifolds, $(\omega_t)_{t \in B}$ is a $C^{\infty}$ family of Hermitian-symplectic metrics on $(X_t)_{t \in B}$. If
\begin{itemize}
\item[(1)] for every $t \in B$ sufficiently close to $0$, $h_{BC,\,t} = h_{BC,\,0}$,
\item[(2)] for every $t \in B$ sufficiently close to $0$, $h_{\bar{\partial},\, t} = h_{\bar{\partial},\, 0} $, 
\item[(3)] for every $t \in B\setminus\{0\}$ and sufficiently close to $0$, the $(0,\,2)$-torsion form  $\rho_{\omega_t} ^{0,\,2}$ is $\bar\partial_t$-exact. 
\end{itemize}
Then 
\begin{itemize}
\item[(a)] the family $(\rho_{\omega_t} ^{0,\,2})_{t  \in B}$ is a $C^{\infty}$ family of $(0,\,2)$-forms,
\item[(b)] for $t=0$, the $(0,\,2)$-torsion form  $\rho_{\omega_0} ^{0,\,2}$ is $\bar\partial_0$-exact. 
\end{itemize}
\end{Prop}
Before giving the proof of Proposition \ref{prop: close torsion dbar exact}, we recall the following lemma which will be used in the proof.
\begin{Lem}\label{dbar-exact sol} (\cite{Pop19}, p.28)
Let $(X, \omega)$ be an n-dimensional compact Hermitian manifold Fix $q \in \lbrace0, \ldots, n\rbrace$.
For every $\rho \in \bar\partial(C^{\infty} _{0,q}(X,T^{1,0}X))$, the minimal $L^2$-norm solution of the equation
\begin{equation}\label{eqn:dbar-harm}
\bar\partial\varphi=\rho
\end{equation}
is given by the following Neumann formula
\begin{equation}\label{eqn:dbar-harm sol}
\varphi= \bar\partial^{\star}(\Delta_{\bar\partial})^{-1}\rho,
\end{equation}
where $(\Delta_{\bar\partial}) ^{-1}$ is the Green operator of the $\bar\partial$-Laplacian $\Delta_{\bar\partial}$ induced by $\omega$. 
\end{Lem}
\textbf{Proof of Proposition \ref{prop: close torsion dbar exact}.} From equation (\ref{eqn:trosion}) one sees that the $\rho_{\omega_t} ^{2,\,0}$ has the following form
\begin{equation}\label{eqn:Neumann_torsion_H-S01}
\rho_{\omega_t} ^{2,\,0}= -\Delta_{BC,\,t}^{-1}[\bar\partial_t ^\star\partial_t \omega_t + \bar\partial_t ^\star\partial_t \partial_t ^\star\partial_t\omega_t],
\end{equation}
for all $t \in B$. By conjugating the above equation one can see that
\begin{equation}\label{eqn:Neumann_conjtorsion_H-S01}
\rho_{\omega_t} ^{0,\,2}= -\Delta_{BC,\,t}^{-1}[\partial_t ^\star\bar\partial_t \omega_t + \partial_t ^\star\bar\partial_t \bar\partial_t ^\star\bar\partial_t\omega_t],
\end{equation}
for all $t \in B$. Note that in (\ref{eqn:Neumann_conjtorsion_H-S01}) we used the fact that $\Delta_{BC}=\overline{\Delta_{BC}}$. Since $h_{BC,\, t} = h_{BC,\, 0 }$ for $t$ sufficiently close to the origin, by Theorem \ref{Kodaira}.(ii), the family $(\Delta_{BC,\,t} ^{-1})_{t \in B}$ of linear operators, acting on $(0,\,2)$-form, is a $C^{\infty}$ family of linear operators. This means that the family $(\rho_{\omega_t} ^{0,\,2})_{t  \in B}$ is a $C^{\infty}$ family of $(0,\,2)$-forms. In particular $\rho_{\omega_t} ^{0,2} \rightarrow \rho_{\omega_0} ^{0,2}$, when $t \longrightarrow 0$. This proves (a).\\
By assumption $(3)$ for every $t \in B\backslash \lbrace 0\rbrace$, $\rho_{\omega_t} ^{0,\,2}$ is $\bar\partial_t$-exact. So after possibly shrinking $B$ about the origin the following equation 
\begin{equation}\label{eqn: e2-exact}
\rho_{\omega_t} ^{0,\,2} =\bar\partial_t \beta_t
\end{equation}
has at least one solution $\beta_t$ in $C^{\infty} _{0,1}(X_t,\mathbb{C})$ for every $t \in B\backslash \lbrace 0\rbrace$. By Lemma \ref{dbar-exact sol}, we are able to choose the unique solution among such $\beta_t$ with the minimal $L^2$-norm induced by $\omega_t$. Hence by equation (\ref{eqn:dbar-harm sol}), the minimal $L^2$-norm solution of equation (\ref{eqn: e2-exact}) has the following form
\begin{equation}\label{beta_t min}
\beta^{\min} _t=\bar\partial^{\star} _t (\Delta_{\bar\partial,\, t})^{-1}\rho_{\omega_t} ^{0,\,2} \stackrel{\text{(I)}}{=} (\Delta_{\bar\partial,\, t})^{-1}\bar\partial^{\star} _t\rho_{\omega_t} ^{0,\,2} .
\end{equation}
Where $(I)$ is implied as follows
\begin{equation*}
\bar\partial^{\star}\Delta_{\bar\partial}=\bar\partial^{\star}(\bar{\partial}\bar\partial^{\star}+\bar\partial^{\star}\bar{\partial})= \bar{\partial}^{\star}\bar\partial\bar{\partial}^{\star}=\bar{\partial}^{\star}\bar{\partial}\bar{\partial}^{\star}  +\bar{\partial}\bar\partial^{\star}\bar{\partial}^{\star}=
(\bar{\partial}\bar\partial^{\star}+\bar\partial^{\star}\bar{\partial})\bar{\partial}^{\star}= \Delta_{\bar\partial}\bar{\partial}^{\star}.
\end{equation*}
By assumption $(2)$ after possibly shrinking $B$ about the origin, $h_{\bar{\partial},\, t} = h_{\bar{\partial},\, 0} $ for every  $t \in B$. Therefore by Theorem \ref{Kodaira} (ii) the family $(\Delta ^{-1} _{\bar\partial,\, t})_{t \in B}$ is a $C^{\infty}$ family of linear operators acting on $(0,\,1)$-forms. On the other hand from (a) one can imply that the family $(\rho_{\omega_t} ^{0,\,2})_{t  \in B}$ is a $C^{\infty}$ family of $(0,\,2)$-forms. Hence there exists a $\beta_0=\Delta_{\bar{\partial},\,0}\bar\partial_0 ^{\star}\rho_{\omega_0} ^{0,\,2} \in  C^{\infty} _{0,\,1}(X_0,\mathbb{C})$ such that the family $(\beta^{\min} _t)_{t \in B}$ is a $C^{\infty}$ family of $(0,\, 1)$-forms. In other words 
\begin{equation}\label{beta0}
\lim_{t \rightarrow 0} \beta^{\min} _t= \beta_0. 
\end{equation}
From equations (\ref{beta_t min}) and (\ref{beta0}) we get
\begin{eqnarray}
\noindent \bar\partial_0 \beta_0 =\bar\partial_0\lim_{t \rightarrow 0} \beta^{\min} _t \stackrel{\text{(I)}}{=} \lim_{t \rightarrow 0} \bar\partial_t \beta^{\min} _t \stackrel{\text{(II)}}{=}  \lim_{t \rightarrow 0} \bar\partial_t \bar\partial^{\star} _t(\Delta_{\bar\partial,\, t})^{-1}\rho_{\omega_t} ^{0,\,2}\stackrel{\text{(III)}}{=}  \lim_{t \rightarrow 0} \rho_{\omega_t} ^{0,\,2}.
\end{eqnarray}
In the above equation, (I) comes from the fact that the family $(\bar\partial _t)_{\ t \in B}$ is a $C^{\infty}$ family of smooth linear operators so it commutes with $\lim$, (II) comes from the definition of $\beta^{\min} _t$ in equation (\ref{beta_t min}) and finally, we have (III) because
\begin{eqnarray}\label{eqn:dbar star dbar commutes delta dol}
\nonumber \rho_{\omega_t} ^{0,\,2} =\Delta_{\bar\partial,\, t}(\Delta_{\bar\partial,\, t})^{-1}\rho_{\omega_t} ^{0,\,2}= (\bar\partial^{\star} _t\bar\partial_t+\bar\partial_t \bar\partial^{\star} _t)(\Delta_{\bar\partial,\, t})^{-1} \rho_{\omega_t} ^{0,\,2}\\= 
\bar\partial_t \bar\partial^{\star} _t(\Delta_{\bar\partial,\, t})^{-1} \rho_{\omega_t} ^{0,\,2} + \bar\partial^{\star} _t\bar\partial_t(\Delta_{\bar\partial,\, t})^{-1} \rho_{\omega_t} ^{0,\,2},
\end{eqnarray} 
first note that 
\begin{equation*}
\bar\partial^{\star}\bar\partial\Delta_{\bar\partial}= \bar\partial^{\star}\bar\partial(\bar\partial^{\star}\bar\partial+ \bar\partial\bar\partial^{\star})= \bar\partial^{\star}\bar\partial\bar\partial^{\star}\bar\partial= \bar\partial^{\star}\bar\partial\bar\partial^{\star}\bar\partial+ \bar\partial\bar\partial^{\star}\bar\partial^{\star}\bar\partial= (\bar\partial^{\star}\bar\partial+ \bar\partial\bar\partial^{\star})\bar\partial^{\star}\bar\partial= \Delta_{\bar\partial}\bar\partial^{\star}\bar\partial,
\end{equation*}
so $(\Delta_{\bar\partial,\, t})^{-1} $ commutes with $\bar\partial^{\star}\bar\partial$, hence in equation (\ref{eqn:dbar star dbar commutes delta dol}) one gets
$$  \bar\partial^{\star} _t\bar\partial_t(\Delta_{\bar\partial,\, t})^{-1} \rho_{\omega_t} ^{0,\,2}= (\Delta_{\bar\partial,\, t})^{-1} \bar\partial^{\star} _t\bar\partial_t \rho_{\omega_t} ^{0,\,2}  $$
and since $\bar\partial_t \rho_{\omega_t} ^{0,\,2}=0$, $ (\Delta_{\bar\partial,\, t})^{-1}\bar\partial^{\star} _t\bar\partial_t \rho_{\omega_t} ^{0,\,2} $ vanishes, so
\begin{equation*}
\bar\partial_t \bar\partial^{\star} _t(\Delta_{\bar\partial,\, t})^{-1} \rho_{\omega_t} ^{0,\,2} + \bar\partial^{\star} _t\bar\partial_t(\Delta_{\bar\partial,\, t})^{-1} \rho_{\omega_t} ^{0,\,2}=\bar\partial_t \bar\partial^{\star} _t(\Delta_{\bar\partial,\, t})^{-1} \rho_{\omega_t} ^{0,\,2}.
\end{equation*}
From (a) one can see that the family $(\rho_{\omega_0} ^{0,\,2})_{t \in B}$ is a $C^{\infty}$ family of $(0,\,2)$-forms. Which means that 
$$\lim_{t \rightarrow 0} \rho_{\omega_t} ^{0,\,2} =\rho_{\omega_0} ^{0,\,2}.$$
This proves (b). \hfill $\Box$\\
In Proposition \ref{prop: close torsion dbar exact}, not only did we prove that the $\bar\partial$-exactness for the family $\rho_{\omega_0} ^{0,\,2}$ is a closed property under holomorphic deformations but also we showed that the family $(\beta^{\text{min}} _t)_{t \in B}$ of minimal $L^2 _{\omega_t}$ solutions is a $C^{\infty}$ family of $(1,\,0)$-forms and the existence of a minimal $L^2 _{\omega_t}$ solution is closed property under holomorphic deformations.\\ 
From now on we focus on the Dinew-Popovici energy functional $F$ defined in section \ref{section:Introduction} and its critical points.
In the following we give a proof for Proposition \ref{prop: close critical F_t}. We show that for a fix Hermitian-symplectic metric $\omega$ being a critical point for Dinew-Popovici energy functional $F$ is a closed property under holomorphic deformations. 
\\ \textbf{Proof of Proposition \ref{prop: close critical F_t}.}
From  Proposition \ref{Prop:F-tilde_properties} for every $t \in B$ and for every $(1,1)$-form $\gamma_t=\bar{\partial_t}u_t+\partial_t\bar u_t$, one gets
\begin{eqnarray}
\nonumber(d_{\omega_t}F_t)(\gamma_t)  = -2\,\mbox{Re}\,\langle\langle u_t,\,\bar\partial_t ^{\star}\omega_t\rangle\rangle_{\omega_t} + 2\,\mbox{Re}\,\int\limits_{X_t} u_t\wedge\rho_{\omega_t} ^{2,\,0}\wedge\overline{\rho_{\omega_t} ^{2,\,0}}\wedge\bar\partial_t \frac{\omega_t ^{n-3}}{(n-3)!}.
 \end{eqnarray}
Since $\omega_t$ is a critical point of $F_t$ for $t \in B\backslash \lbrace 0\rbrace$, $(d_{\omega_t}F_t)(\gamma_t)=0$ for every $\gamma_t \in C^{\infty}_{1,\,1}(X_t,\mathbb{C})$. Also by assumption $(2)$ and Proposition \ref{prop: close torsion dbar exact} the family $ (\rho_{\omega_t} ^{0,\,2})_{t \in B}$ is a $C^{\infty}$ family of $(0, \, 2)$-forms. It is obvious that the smooth $J_{t}$-$(1,0)$ forms $u_t$ on $X_t$ determines $d_{\omega_t}F_t$. Define for every $t \in B$,
\begin{eqnarray}\label{eqn:critical-close}
& \nonumber T_t: C^{\infty} _{1, \, 0} (X_t, \mathbb{C}) \longrightarrow \mathbb{R} \hspace{2cm} T_t(u_t)= G_t(u_t)+ H_t(u_t),
\end{eqnarray}
where
\begin{equation}
G_t( u_t)= -2\,\mbox{Re}\,\langle\langle u_t,\,\bar\partial_t ^{\star}\omega_t\rangle\rangle_{\omega_t}
\end{equation}
and 
\begin{equation}
H_t (u_t)= 2 \,\mbox{Re}\,\int\limits_{X_t} u_t\wedge\rho_{\omega_t} ^{2,\,0}\wedge\overline{\rho_{\omega_t} ^{2,\,0}}\wedge\bar\partial_t \frac{\omega_t ^{n-3}}{(n-3)!}.
\end{equation}
In order to prove that $\omega_0$ is a critical point for $F_0$, it is sufficient to show that the family $(T_t)_{t \in B}$ is a $C^{\infty}$ family of linear operators. Therefore it is sufficient to consider a $C^{\infty}$ family of $(1,\,0)$-forms $(u_t)_{t \in B}$ and show that $T_t (u_t)=0$ for all $t \in B$. \\ 
Now suppose that the  $C^{\infty}$ family  of $(1, \, 0)$-forms $(u_t)_{t \in B}$  is given and $B$ is sufficiently shrunk about the origin. 
For every $t \in B$, 
$$\bar{\partial}_t ^{\star}: C^{\infty}_{1,\,1}(X_t,\,\mathbb{C}) \longrightarrow C^{\infty}_{1,\,0}(X_t,\,\mathbb{C}),$$
is a smooth linear operator and the family $(\bar{\partial}_t ^{\star})_{t \in B}$ is a $C^{\infty}$ family of linear operators. Also for every $t \in B$, the map
$$ \langle\langle ~ \, ,\, \omega_t \rangle\rangle_{\omega_t} : C^{\infty}_{1,\,1}(X_t,\,\mathbb{C}) \longrightarrow \mathbb{C}, \hspace{0.75cm}  \langle\langle ~ \, ,\, \omega_t \rangle\rangle_{\omega_t} (\alpha)= \langle\langle \alpha\, ,\, \omega_t \rangle\rangle_{\omega_t}$$
is a smooth linear map and the family $( \langle\langle ~\, ,\, \omega_t \rangle\rangle_{\omega_t})_{t \in B}$ is a $C^{\infty}$ family of linear operators. So for every $t \in B$, $G_t$ is a smooth linear operator and the family $(G_t)_{t \in B}$ is a $C^{\infty}$ family linear operators. In other words
\begin{equation}\label{G_0}
\lim_{t \to 0} G_t(u_t)= G_0(u_0)=-2\,\mbox{Re}\,\langle\langle u_0,\,\bar\partial_0 ^{\star}\omega_0\rangle\rangle_{\omega_0}. 
\end{equation}
We show that the family $(H_t)_{t \in B}$ is a $C^{\infty}$ family of linear operators. First it is obvious that for every $t \in B$, the map 
$$\bar\partial_t: C^{\infty}_{n-3,\,n-3}(X_t,\,\mathbb{C}) \longrightarrow C^{\infty} _{{n-3,\,n-2}}(X_t,\,\mathbb{C}),$$
is a smooth linear operator and the family $(\bar\partial_t)_{t \in B}$ is a $C^{\infty}$ family of linear operators. On the other hand, the  family $(\omega_t)_{t \in B}$  is a $C^{\infty}$ family of metrics, henceforth the family $(\frac{\bar\partial\omega_t ^{n-3}}{(n-3)!})_{t \in B}$ is a $C^{\infty}$ family of $(n-3,\,n-2)$-forms. Also, assumption (2) allows us to employ Proposition \ref{prop: close torsion dbar exact} and say that both families $(\rho_{\omega_t} ^{0,\,2})_{t \in B}$ and $(\overline{\rho_{\omega_t} ^{0,\,2}})_{t \in B}$ are $C^{\infty}$ family of $(2 ,\,0)$-forms and $(0 ,\,2)$-forms respectively. Therefore for every $t \in B$ the map $H_t$ is a smooth real-valued linear map and the family $(H_t)_{t \in B}$ is a  $C^{\infty}$ family of linear operators. In other words,
\begin{equation}\label{H_0}
\lim_{t \to 0} H_t(u_t)= H_0(u_0)=\,\mbox{Re}\,\int\limits_{X_0} u_0 \wedge\rho_{\omega_0} ^{2,\,0}\wedge\overline{\rho_{\omega_0} ^{2,\,0}}\wedge\bar\partial_0 \frac{\omega_0 ^{n-3}}{(n-3)!}.
\end{equation}
The smoothness of $T_t$ for every $t \in B$ is implied by the smoothness of $G_t$ and $H_t$, and by equations (\ref{G_0}) and (\ref{H_0}) one can get
\begin{equation}
T_0(u_0)= \lim_{t \in B} T_t(u_t)= \lim_{t \in B} G_t(u_t)+ \lim_{t \in B} H_t(u_t)=0.
\end{equation}
 Which means that $\omega_0$ is a critical point of $F_0.$ \hfill $\Box$\\
In section \ref{section: back energy functioal} we saw that in dimension 3 the explicit formula for differential of the Dinew-Popovici energy functional $F$ a $\omega$ is simpler than in higher dimensions. For next result of this article we give a proof to Proposition  \ref{intro :prop: critial d-exact}, where we compute the differential of $F$ at $\omega$, when $\omega$ is a fixed Hermitian-symplectic metric on compact complex manifold $X$ of dimension $n$ and the $(2,\,0)$-torsion form $\rho^{2,0} _{\omega}$ is $\partial$-exact.

\textbf{Proof of Proposition \ref{intro :prop: critial d-exact}.} First note that since $\omega \in {\cal S}_{\{\omega_0\}}$, it is a Hermitian-symplectic metric on $X$ so the  $(2,\,0)$-torsion form $\rho_\omega^{2,\,0}$  satisfies, $\bar\partial\omega= -\partial\overline{\rho_\omega^{2,\,0}}$ and $\partial \omega=-\bar\partial \rho_\omega^{2,\,0} $ and $\bar\partial\overline{\rho_\omega^{2,\,0}}=0$. On the other hand, since $X$ is a compact complex manifold it has no boundary so for every $(n-1,\,n)$-form $\alpha$ and every $(n,\,n-1)$-form $\beta$  
$$\int\limits_X \partial \alpha =0  ~~~~\text{and}~~~~ \int\limits_X \bar\partial\beta=0$$
by Stokes's theorem. From (\ref{eqn:differential_F}), one can observe that when $\gamma= \bar\partial\xi+\partial\bar\xi$ the differential at $\omega$ of $F$ evaluated on $\gamma$ is
\begin{equation} \label{dif of f xi}
(d_{\omega}F)(\gamma) = -2\,\mbox{Re}\,\langle\langle \xi,\,\bar\partial^{\star}\omega\rangle\rangle_\omega + 2\,\mbox{Re}\,\int\limits_X \xi\wedge\rho_\omega^{2,\,0}\wedge\overline{\rho_\omega^{2,\,0}}\wedge\bar\partial\omega_{n-3}.
\end{equation}
First we compute $\langle\langle \xi,\,\bar\partial^{\star}\omega\rangle\rangle_\omega$. By the definition of  the $L^2 _{\omega}$ inner product ($u \wedge \star \bar v = \langle u,\, v\rangle_{\omega} dV_{\omega} $), we have
\begin{equation}\label{xi  and omega}
\langle\langle \xi,\,\bar\partial^{\star}\omega\rangle\rangle_\omega = \int\limits_X \langle \xi,\, \bar\partial^{\star} \omega \rangle_{\omega} dV_{\omega}= \int\limits_X \xi \wedge \star\partial^{\star}\omega.
\end{equation}
By standard computation for the Hodge star operator $-\star\star= id$ on odd-degree forms and by equation (\ref{star op on omega }), one gets
\begin{equation}\label{hodge omega}
\star\partial^{\star}\omega= -\star \star\bar\partial \star\omega= \bar\partial\omega_{n-1},
\end{equation}
hence by equations (\ref{xi  and omega}) and (\ref{hodge omega}), 
\begin{equation}\label{xi and omega 2}
\langle\langle \xi,\,\bar\partial^{\star}\omega\rangle\rangle_\omega=  \int\limits_X  \xi \wedge \bar\partial\omega_{n-1}.
\end{equation}
Now equation (\ref{xi and omega 2}) allows us to compute $\langle\langle \xi,\,\bar\partial^{\star}\omega\rangle\rangle_\omega$. We get 
\begin{equation}\label{fist part}
\langle\langle \xi,\,\bar\partial^{\star}\omega\rangle\rangle_\omega=  \int\limits_X  \xi \wedge \bar\partial\omega_{n-1}=  \int\limits_X  \xi \wedge \bar\partial\omega \wedge \omega_{n-2}= -\int\limits_X  \xi \wedge \partial\overline{\rho_\omega^{2,\,0}} \wedge \omega_{n-2}.
\end{equation}
By Stokes's theorem $0=\int\limits_X \partial (\xi \wedge \overline{\rho_\omega^{2,\,0}} \wedge \omega_{n-2})$, so
\begin{equation}
0=\int\limits_X  \partial\xi \wedge \overline{\rho_\omega^{2,\,0}} \wedge \omega_{n-2} -\int\limits_X  \xi \wedge  \partial\overline{\rho_\omega^{2,\,0}} \wedge \omega_{n-2}- \int\limits_X  \xi \wedge \overline{\rho_\omega^{2,\,0}} \wedge \partial \omega_{n-2}.
\end{equation}
Therefore, 
\begin{equation}
 -\int\limits_X  \xi \wedge \partial\overline{\rho_\omega^{2,\,0}} \wedge \omega_{n-2}= \int\limits_X  \xi \wedge \overline{\rho_\omega^{2,\,0}} \wedge \partial \omega_{n-2} - \int\limits_X  \partial\xi \wedge \overline{\rho_\omega^{2,\,0}} \wedge \omega_{n-2}.
\end{equation}
By assumption $ \rho_{\omega}^{2,0}= \partial\xi$ so, 
$$\int\limits_X  \partial\xi \wedge \overline{\rho_\omega^{2,\,0}} \wedge \omega_{n-2}=\int\limits_X  \rho_\omega^{2,\,0} \wedge \overline{\rho_\omega^{2,\,0}} \wedge \omega_{n-2}.$$
Since $\rho_{\omega}^{2,0}$ is a primitive form of bidegree $(2,\,0)$, we can apply (\ref{eqn:primi formula for fomrs}) and we get:
$$\int\limits_X  \rho_\omega^{2,\,0} \wedge \overline{\rho_\omega^{2,\,0}} \wedge \omega_{n-2}= \int\limits_X \rho_\omega^{2,\,0} \wedge  \star \overline{\rho_\omega^{2,\,0}}= \langle\langle \rho_\omega^{2,\,0},\,\rho_\omega^{2,\,0}\rangle\rangle_\omega= \Vert \rho_{\omega}^{2,0}\Vert^2 _{\omega}. $$ 
Hence 
\begin{equation}\label{sec stoke}
 -\int\limits_X  \xi \wedge \partial\overline{\rho_\omega^{2,\,0}} \wedge \omega_{n-2}= \int\limits_X  \xi \wedge \overline{\rho_\omega^{2,\,0}} \wedge \partial\omega_{n-2} - \Vert \rho_{\omega}^{2,0}\Vert^2 _{\omega}.
\end{equation}
Now, the goal is to compute $ \int\limits_X  \xi \wedge \overline{\rho_\omega^{2,\,0}} \wedge \partial\omega_{n-2} $ in equation (\ref{sec stoke}). We get:
\begin{equation} \label{sec stoke3}
 \int\limits_X  \xi \wedge \overline{\rho_\omega^{2,\,0}} \wedge \partial\omega_{n-2}=  \int\limits_X  \xi \wedge \overline{\rho_\omega^{2,\,0}} \wedge \partial\omega\wedge\omega_{n-3}= - \int\limits_X  \xi \wedge \overline{\rho_\omega^{2,\,0}} \wedge \bar\partial \rho_\omega^{2,\,0}\wedge\omega_{n-3}.
\end{equation}
Again, by Stokes's theorem, $ \int\limits_X \bar\partial (\xi \wedge \overline{\rho_\omega^{2,\,0}} \wedge \rho_\omega^{2,\,0} \wedge \omega_{n-3})=0 $ and because $\bar\partial  \overline{\rho_\omega^{2,\,0}}=0 $, we have 
\begin{equation}\label{sec stoke1}
0= \int\limits_X \bar\partial\xi \wedge \overline{\rho_\omega^{2,\,0}} \wedge \rho_\omega^{2,\,0} \wedge \omega_{n-3}-
\int\limits_X \xi \wedge \overline{\rho_\omega^{2,\,0}} \wedge \bar\partial\rho_\omega^{2,\,0} \wedge \omega_{n-3}-
\int\limits_X \xi \wedge \overline{\rho_\omega^{2,\,0}} \wedge \rho_\omega^{2,\,0} \wedge \bar\partial\omega_{n-3}.
\end{equation}
Therefore, from (\ref{sec stoke1}) and (\ref{sec stoke3})  one can deduce the following equation 
\begin{equation}\label{last eqn}
\int\limits_X  \xi \wedge \overline{\rho_\omega^{2,\,0}} \wedge \partial \omega_{n-2}= -\int\limits_X \bar\partial\xi \wedge \overline{\rho_\omega^{2,\,0}} \wedge \rho_\omega^{2,\,0} \wedge \omega_{n-3}+
\int\limits_X \xi \wedge \overline{\rho_\omega^{2,\,0}} \wedge \rho_\omega^{2,\,0} \wedge \bar\partial\omega_{n-3}.
\end{equation}
By putting equations (\ref{fist part}), (\ref{sec stoke}) and (\ref{last eqn}) together we see that
\begin{equation}\label{inner eqn}
-2\,\mbox{Re}\langle\langle \xi,\,\bar\partial^{\star}\omega\rangle\rangle_\omega= 2\,\mbox{Re}\int\limits_X \bar\partial\xi \wedge \overline{\rho_\omega^{2,\,0}} \wedge \rho_\omega^{2,\,0} \wedge \omega_{n-3}-2\,\mbox{Re}
\int\limits_X \xi \wedge \overline{\rho_\omega^{2,\,0}} \wedge \rho_\omega^{2,\,0} \wedge \bar\partial\omega_{n-3}+2  \Vert \rho_{\omega}^{2,0}\Vert^2 _{\omega}.
\end{equation}
By adding $2\,\mbox{Re}\,\int\limits_X \xi \wedge\rho_\omega^{2,\,0}\wedge\overline{\rho_\omega^{2,\,0}}\wedge\bar\partial\omega_{n-3}$ to equation (\ref{inner eqn}) and by using (\ref{eqn:differential_F}) with $u =\xi$ we get the formula (\ref{der F wrt xi}). This proves the proposition. \hfill $\Box$\\

In formula (\ref{der F wrt xi}), $2\mbox{Re} \int\limits_X  \bar\partial \xi \wedge \rho^{2,0} \wedge  \overline{\rho_\omega^{2,\,0}} \wedge \omega_{n-3}$ is signless in general. However, if it supposes to be non-negative one sees immediately that $\omega$ is a K\"ahler metric whenever it is a critical point for $F$. In the following proof, we show that if $\bar{\partial}\xi$ is a weakly semi-positive $(1,\,1)$-form then $2\mbox{Re} \int\limits_X  \bar\partial \xi \wedge \rho^{2,0} \wedge  \overline{\rho_\omega^{2,\,0}} \wedge \omega_{n-3}$ is non-negative.

\textbf{Proof of Corollary \ref{intro :prop: critial d-exact}.} Since positivity is a pointwise property, one can fix a point $x \in X$  and local coordinates $(z_1,\ldots,z_n)$ centered at $x$ such that $\omega$ has the following shape
$$\omega =  \sum idz_i \wedge d\bar{z}_i ~~~~~~~~~ \text{at}~x$$
In particular, $\omega$ is a strongly strictly positive $(1,\,1)$-form. By Definition \ref{positivty} (1) and Proposition \ref{prod of sp}, $\omega_{n-3}$ is a strongly strictly positive $(n-3,\,n-3)$-form. On the other hand by Example 1.2 of \cite{Dem}, for every $p \in \{1,\dots,n\}$ and any $(p,\,0)$-form $\beta$, the $(p,\,p)$-form $i^{p^2}\beta \wedge \bar\beta$ is weakly strictly positive. Hence  the $(2,\,2)$-form 
$$i^{4} \rho_{\omega} ^{2,0} \wedge \overline{\rho_{\omega} ^{2,0}}= \rho_{\omega} ^{2,0} \wedge \overline{\rho_{\omega} ^{2,0}}$$ 
is weakly strictly positive. \\
Since $\bar\partial\xi$ is weakly semi-positive $(1,\,1)$-form, there exist real non negative functions  $c_1,\dots,\,c_n $ and $(1,\,0)$-forms $\alpha_1,\dots,\, \alpha_n$ such that 
$$ \bar\partial\xi = \sum c_k i \alpha_k \wedge \bar{\alpha}_k.$$
Therefore 
\begin{eqnarray}
\nonumber \bar\partial \xi \wedge \rho^{2,0} \wedge  \overline{\rho_\omega^{2,\,0}} \wedge \omega_{n-3}= 
 \sum c_k i \alpha_k \wedge \bar{\alpha}_k \wedge \rho^{2,0} \wedge  \overline{\rho_\omega^{2,\,0}} \wedge \omega_{n-3}\\
 =  \sum c_k \rho^{2,0} \wedge  \overline{\rho_\omega^{2,\,0}} \wedge i\alpha_k \wedge \bar{\alpha}_k \wedge \omega_{n-3}.
\end{eqnarray}  
Note that by Definition \ref{positivty} $\alpha_k \wedge \bar{\alpha}_k $ is strongly strictly positive $(1,\,1)$-form for all $k \in \{ 1,\dots,\,n\}$. By Definition \ref{positivty} (2),  $c_k \rho^{2,0} \wedge  \overline{\rho_\omega^{2,\,0}} \wedge i\alpha_k \wedge \bar{\alpha}_k \wedge \omega_{n-3}$ is a weakly semi-positive $(n,\,n)$-form. Hence 
$$2\mbox{Re} \int\limits_X \sum c_k \rho^{2,0} \wedge  \overline{\rho_\omega^{2,\,0}} \wedge i\alpha_k \wedge \bar{\alpha}_k \wedge \omega_{n-3}
= \sum 2\mbox{Re} \int\limits_Xc_k \rho^{2,0} \wedge  \overline{\rho_\omega^{2,\,0}} \wedge i\alpha_k \wedge \bar{\alpha}_k \wedge \omega_{n-3} \geqslant 0$$
This proves the Corollary. \hfill $\Box$\\

\vspace{3ex}

\noindent Institut de Math\'ematiques de Toulouse, Universit\'e Paul Sabatier,

\noindent 118 route de Narbonne, 31062 Toulouse, France

\noindent Email: Soheil.Erfan@math.univ-toulouse.fr

\end{document}